\documentclass[11pt,reqno]{amsart}
\usepackage{amsmath,amssymb,geometry,color,tikz-cd}
\usepackage[backref,pagebackref,pdftex,hyperindex]{hyperref}
\numberwithin{equation}{section}%
\usepackage[alphabetic]{amsrefs}
\usepackage{amssymb}
\usepackage{bbm}
\usepackage{enumitem}
\usepackage{upgreek}

\newtheorem{proposition}{Proposition}[section]
\newtheorem{theorem}[proposition]{Theorem}
\newtheorem{lemma}[proposition]{Lemma}
\newtheorem{corollary}[proposition]{Corollary}

\theoremstyle{definition}
\newtheorem{remark}[proposition]{Remark}
\newtheorem{definition}[proposition]{Definition}

\newtheorem{claim}[proposition]{Claim}

\title{Weakly special test configurations of log canonical Fano varieties}
\author{Guodu Chen and Chuyu Zhou}

\address{Institute for Theoretical Sciences, Westlake Institute for Advanced Study, Westlake University, Hangzhou, Zhejiang, 310024, China}
\email{chenguodu@westlake.edu.cn}
	
\address{\'Ecole Polytechnique F\'ed\'erale de Lausanne (EPFL), MA C3 615, Station 8, 1015 Lausanne, Switzerland}
\email{chuyu.zhou@epfl.ch}
\date{} 


\newcommand{\Fut}{{\rm{Fut}}}

\newcommand{\ord}{{\rm {ord}}}
\newcommand{\tc}{{\rm {tc}}}
\newcommand{\vol}{{\rm {vol}}}
\newcommand{\lct}{{\rm {lct}}}

\newcommand{\Spec}{{\rm {Spec}}}
\newcommand{\Proj}{{\rm{Proj}}}

\newcommand{\dt}{{\rm {dt}}}

\newcommand{\Ding}{{\rm {Ding}}}
\newcommand{\NA}{{\rm {NA}}}
\newcommand{\Image}{{\rm {Image}}}

\usepackage{todonotes}

\newcommand{\bA}{\mathbb{A}}

\newcommand{\bC}{\mathbb{C}}

\newcommand{\bN}{\mathbb{N}}

\newcommand{\bP}{\mathbb{P}}
\newcommand{\bQ}{\mathbb{Q}}
\newcommand{\bR}{\mathbb{R}}

\newcommand{\bZ}{\mathbb{Z}}

\newcommand{\mD}{\mathcal{D}}
\newcommand{\mE}{\mathcal{E}}
\newcommand{\mF}{\mathcal{F}}

\newcommand{\mI}{\mathcal{I}}

\newcommand{\mL}{\mathcal{L}}

\newcommand{\mO}{\mathcal{O}}

\newcommand{\mW}{\mathcal{W}}
\newcommand{\mX}{\mathcal{X}}
\newcommand{\mY}{\mathcal{Y}}
\newcommand{\mZ}{\mathcal{Z}}

\newcommand{\fD}{\mathbf{D}}
\newcommand{\fE}{\mathbf{E}}

\newcommand{\fL}{\mathbf{L}}

\newcommand{\tX}{\tilde{X}}

\newcommand{\tZ}{\tilde{Z}}

\newcommand{\ka}{\mathfrak{a}}

\begin{document}

\begin{abstract}
Let $(X,\Delta)$ be a strictly lc log Fano pair, we show that every lc place of complements of $(X,\Delta)$ is dreamy, and there exists a correspondence between weakly special test configurations of $(X,\Delta;-K_X-\Delta)$ and lc places of complements of $(X,\Delta)$.
\end{abstract}

\maketitle
\tableofcontents

\section{Introduction}

Throughout the paper we work over $\bC$. We say that $(X,\Delta)$ is a log pair if $X$ is a projective normal variety and $\Delta$ is an effective $\bQ$-divisor on $X$ such that $K_X+\Delta$ is $\bQ$-Cartier. Let $(X,\Delta)$ be a log pair such that $-K_X-\Delta$ is ample, we say that $(X,\Delta)$ is an \emph{lc} (resp. \emph{klt}) \emph{log Fano} pair if it admits log canonical (resp. Kawamata log terminal) singularities.

\begin{definition}
Let $(X,\Delta)$ be an lc log Fano pair. We say that a test configuration (see Definition \ref{defn:tc}) $(\mX,\Delta_\tc; \mL:=-K_{\mX}-\Delta_\tc)\to \bA^1$ of $(X,\Delta; -K_X-\Delta)$ is \emph{weakly special} if  $\mX_0$ is irreducible, and $(\mX, \mX_0+ \Delta_{\tc,0})$ is log canonical.
\end{definition}

Note that our definition of weakly special test configuration is different from the usual one (see \cite[Appendix A]{BLX19}) as here we require the centeral fiber to be irreducible.

Let $(X,\Delta)$ be an lc log Fano pair. We say that $E$ is a prime divisor over $X$ if there is a projective normal birational model $f: Y\to X$ such that $E$ is a prime divisor on $Y$.
We say that a prime divisor $E$ over $X$ is an \emph{lc place of complements of $(X,\Delta)$} if there is an effective $\bQ$-divisor $D\sim_\bQ -K_X-\Delta$ such that $(X,\Delta+D)$ is log canonical and $E$ is an lc place of $(X,\Delta+D)$. The main purpose of this article is to show the following correspondence between weakly special test configurations of $(X,\Delta; -K_X-\Delta)$ and lc places of complements of $(X,\Delta)$.

\begin{theorem}\label{thm: main1}
Let $(X,\Delta)$ be a strictly lc log Fano pair. Suppose $E$ is an lc place of complements, then $E$ induces a non-trivial weakly special test configuration  $(\mX, \Delta_\tc; -K_{\mX}-\Delta_\tc)$ of $(X,\Delta; -K_X-\Delta)$ such that the restriction of $\ord_{\mX_0}$ to the function field $K(X)$ is $c\cdot\ord_E$ for some $c\in \bZ^+$; conversely, if $(\mX,\Delta_\tc; -K_{\mX}-\Delta_\tc)$ is a non-trivial weakly special test configuration of $(X,\Delta;-K_X-\Delta)$, then the restriction of $\ord_{\mX_0}$ to $K(X)$ is an lc place of complements.
\end{theorem}

\begin{remark}
For a klt log Fano pair $(X,\Delta)$, the work \cite[Appendix A]{BLX19} gives the correspondence between weakly special test configurations of $(X,\Delta; -K_X-\Delta)$ and lc places of complements of $(X,\Delta)$. The above result is a generalization of this correspondence to lc setting.
\end{remark}

To show such a correspondence, some finite generation properties play a key role. We need the following definition of a dreamy divisor over an lc log Fano pair.

\begin{definition}{(\cite[Definition 1.3]{Fuj19})}
Let $(X,\Delta)$ be an lc log Fano pair and $E$ a prime divisor over $X$. We say that $E$ is a \emph{dreamy} divisor (or $E$ is dreamy) if the following $\bN^2$-graded algebra 
$$\oplus_{k\in \bN}\oplus_{j\in \bN} H^0(Y, -krf^*(K_X+\Delta)-jE),$$
is finitely generated for some $r\in \bZ^+$ such that $-r(K_X+\Delta)$ is Cartier, where $f:Y\to X$ is a projective normal birational model such that $E$ is a prime divisor on $Y$.
\end{definition}

The second main result is the following dreamy property of lc places of complements for a strictly lc log Fano pair.

\begin{theorem}\label{thm: main2}
Let $(X,\Delta)$ be a strictly lc log Fano pair, then every lc place of complements is dreamy.
\end{theorem}

\begin{remark}
We only consider strictly lc case since the result is well-known for a klt log Fano pair. 
\end{remark}

\noindent
{\bf Acknowledgements}:
The authors would like to thank Ziquan Zhuang for valuable discussions on Theorem \ref{key lemma}. G. Chen is supported by the China post-doctoral grants BX2021269 and 2021M702925. C. Zhou is supported by grant European Research Council (ERC-804334). Finally, the authors are grateful to the referees for many valuable comments and suggestions.

\section{Preliminaries}

In this section, we give a brief introduction to the concepts of generalized Futaki invariant, Ding invariant, and log canonical slope. The readers may refer to \cite{Fuj19, XZ20b,Li19} for more details on these concepts. For various types of singularities in birational geometry such as klt, lc, dlt, etc. we refer the readers to \cite{KM98, Kollar13}.

\subsection{Test configuration}

\begin{definition}\label{defn:tc}
Let $(X,\Delta)$ be a log pair of dimension $d$ and $L$ an ample $\bQ$-line bundle on $X$. A \emph{test configuration} $\pi: (\mX,\Delta_\tc;\mL)\to \bA^1$ is a family over $\bA^1$ consisting of the following data:
\begin{enumerate}
\item $\pi: \mX\to \bA^1$ is a projective flat morphism from a normal variety $\mX$, $\Delta_\tc$ is an effective $\bQ$-divisor on $\mX$ and $\mL$ is a relatively ample $\bQ$-line bundle on $\mX$,
\item the family $\pi$ admits a $\bC^*$-action which lifts the natural $\bC^*$-action on $\bA^1$ such that $(\mX,\Delta_\tc; \mL)\times_{\bA^1}\bC^*$ is $\bC^*$-equivariantly isomorphic to $(X, \Delta; L)\times_{\bA^1}\bC^*$.
\end{enumerate}
\end{definition}

We denote $\left(\bar{\mX}, \bar{\Delta}_\tc;\bar{\mL}\right)\to \bP^1$ to be the natural compactification of the original test configuration, which is obtained by glueing $(\mX, \Delta_\tc;\mL)$ and $(X,\Delta;L)\times \left(\bP^1\setminus 0\right)$ along their common open subset $(X, \Delta;L)\times \bC^*$.

\begin{definition}\label{tc}
The \emph{generalized Futaki invariant} of a test configuration $(\mX,\Delta_\tc;\mL)$ is defined as follows:
$$\Fut(\mX,\Delta_\tc;\mL):=\frac{\left(K_{\bar{\mX}/\bP^1}+\bar{\Delta}_\tc\right)\cdot\bar{\mL}^d}{L^d} -\frac{d}{d+1}\cdot\frac{\left(\left(K_X+\Delta\right)L^{d-1}\right)\bar{\mL}^{d+1}}{\left(L^d\right)^2}.$$
In particular, if $(X,\Delta;-K_X-\Delta)$ is an lc log Fano pair of dimension $d$, then the generalized Futaki invariant can be expressed as
$$\Fut(\mX,\Delta_\tc;\mL)=\frac{d\bar{\mL}^{d+1}}{(d+1)(-K_X-\Delta)^d} +\frac{\bar{\mL}^d\cdot(K_{\bar{\mX}/\bP^1}+\bar{\Delta}_\tc)}{(-K_X-\Delta)^d}.$$
\end{definition}

\begin{definition}
Let $(X,\Delta)$ be an lc log Fano pair, the \emph{Ding invariant} of a test configuration $(\mX,\Delta_\tc;\mL)$ of $(X,\Delta;-K_X-\Delta)$ is defined as follows
$$\Ding(\mX,\Delta_\tc;\mL):=-\frac{\bar{\mL}^{d+1}}{(d+1)(-K_X-\Delta)^d} +\lct\left(\mX, \Delta_\tc+\mD_{(\mX,\Delta_\tc;\mL)};\mX_0\right)-1, $$
where $\mD_{(\mX,\Delta_\tc;\mL)}\sim_\bQ -(K_{\mX}+\Delta_\tc+\mL)$ is a $\bQ$-divisor on $\mX$ whose support is contained in $\mX_0$ (see \cite[Definition 3.1]{Fuj18}). As Ding invariant does not change if we replace $\mL$ (resp. $\mD_{(\mX,\Delta_\tc;\mL)}$) with $\mL+a\mX_0$ (resp. $\mD_{(\mX,\Delta_\tc;\mL)}-a\mX_0$), sometimes we just write $$\mD_{(\mX,\Delta_\tc;\mL)}= -(K_{\mX}+\Delta_\tc+\mL) .$$
\end{definition}

\begin{remark}\label{coincidence}
For an lc log Fano pair $(X,\Delta)$, it is not hard to see that the generalized Futaki invariant coincides with Ding invariant for a weakly special test configuration.
\end{remark}

Let $(X,\Delta)$ be an lc log Fano pair of dimension $d$, we define the following two well-known invariants, i.e., \emph{log discrepancy} and \emph{$S$-invariant} of a prime divisor $E$ over $X$:
$$A_{X,\Delta}(E):=\ord_{E}(K_Y-f^*(K_X+\Delta))+1, \text{ and }$$
$$S_{X,\Delta}(E):=\frac{1}{(-K_X-\Delta)^d}\int_0^\infty \vol(-f^*(K_X+\Delta)-tE)\dt, $$
where $f: Y\to X$ is a projective normal birational model such that $E$ is a prime divisor on $Y$. More generally, if $L$ is an ample $\bQ$-line bundle on $X$, one can also define the following $S$-invariant associated to $L$:
$$S_{X,L}(E):=\frac{1}{L^d}\int_0^\infty \vol(f^*L-tE)\dt.$$

In the work \cite{Fuj19}, K. Fujita gives a computation for dreamy test configurations (i.e., test configurations with integral central fibers) of a klt log Fano pair. 

\begin{theorem}{\rm {(\cite[Theorem 5.1]{Fuj19})}}\label{futaki}
Let $(X,\Delta)$ be a klt log Fano pair of dimension $d$ and $(\mX, \Delta_\tc; -K_{\mX}-\Delta_\tc)$ a non-trivial test configuration of $(X,\Delta;-K_X-\Delta)$ with integral central fiber. Then we have the following reformulation of the generalized Futaki inavriant:
$$\Fut(\mX,\Delta_\tc; -K_{\mX}-\Delta_\tc)=c\cdot (A_{X,\Delta}(E)-S_{X,\Delta}(E)), $$
where ${\ord_{\mX_0}}|_{K(X)}=c\cdot \ord_E$ for some prime divisor $E$ over $X$ and some $c\in \bQ^+$.
\end{theorem}

The next result is just a restatement of the above theorem in lc setting. In fact, the proof in \cite{Fuj19} does not really depend on kltness, however, we write it as a corollary of a generalization in Section \ref{sec: 3} (see Theorem \ref{generalization}), as we think the generalization may have independent interest.

\begin{corollary}\label{lc computation}
Let $(X,\Delta)$ be an lc log Fano pair of dimension $d$ and $(\mX, \Delta_\tc; -K_{\mX}-\Delta_\tc)$ a non-trivial test configuration of $(X,\Delta;-K_X-\Delta)$ with integral central fiber. Then we have the following reformulation of the generalized Futaki inavriant:
$$\Fut(\mX,\Delta_\tc; -K_{\mX}-\Delta_\tc)=c\cdot \left(A_{X,\Delta}\left(E\right)-S_{X,\Delta}\left(E\right)\right), $$
where ${\ord_{\mX_0}}|_{K(X)}=c\cdot \ord_E$ for some prime divisor $E$ over $X$ and some $c\in \bQ^+$.
\end{corollary}

\begin{proof}
It is implied by Theorem \ref{generalization} by taking $L=-K_X-\Delta$.
\end{proof}

\subsection{Log canonical slope}

Let $(X,\Delta)$ be a log pair and $L$ an ample line bundle on $X$. Denote by $R:=\oplus_{k\in \bN} R_k$ and $R_k:= H^0(X,kL)$. We consider the filtration on $R$ as in the following definition, which is also called linearly bounded multiplicative filtration, see \cite{Fuj19, XZ20b}.

\begin{definition}\label{def: filtration}
By a \emph{filtration} $\mF$ of $R$, we mean the data of a family of $\mathbb{C}$-vector subspaces
$$\mF^\lambda R_k\subset R_k$$
for $k\in \bN$ and $\lambda\in \bR$, satisfying:
\begin{enumerate}
\item (decreasing) $\mF^\lambda R_k\subset \mF^{\lambda'}R_k$ if $\lambda\geq \lambda'$;
\item (left continuous) $\mF^\lambda R_k=\cap_{\lambda'<\lambda}\mF^{\lambda'}R_k$ for all $\lambda\in \bR$;
\item (linearly bounded) there exist $e_-, e_+\in \bR$ such that $\mF^{kx}R_k=0$ for all $x\geq e_+$ and $\mF^{kx}R_k=R_k$ for all $x\leq e_-$;
\item (multiplicative) $\mF^\lambda R_k\cdot \mF^{\lambda'}R_{k'}\subset \mF^{\lambda+\lambda'}R_{k+k'}$.
\end{enumerate}
\end{definition}

A filtration $\mF$ of $R$ is called a \emph{$\bZ$-filtration} if $\mF^\lambda R_k=\mF^{\ulcorner \lambda \urcorner}R_k$ for all $k\in \bN$ and $\lambda\in \bR$. It is called an \emph{$\bN$-filtration} if in addition $\mF^0R_k=R_k$ for all $k\in \bN$. For example, given a valuation $v: K(X)^*\to \bR$ on $X$, it induces a filtration of $R$ by setting $\mF^\lambda R_k:=\{s\in R_k\mid\textit{$v(s)\geq \lambda$}\}$. If $v=c\cdot\ord_E$ for some $c\in \bZ^+$ and some prime divisor $E$ over $X$, then it induces an $\bN$-filtration.

\begin{definition}
Let $\mF$ be a filtration of $R$. For $k\in \bN$ and $\lambda\in \bR$ we set
$$I_{k,\lambda}= I_{k,\lambda}(\mF):={\rm Im}\left(\mF^\lambda R_k\otimes \mO_X\left(-kL\right)\to \mO_X \right),$$
where the map is naturally induced by the evaluation $H^0(X,kL)\otimes \mO_X(-kL)\to \mO_X$. Define $I_k^{(t)}=I_k^{(t)}(\mF):=I_{k,kt}(\mF)$ for $t\in \bR$, then $I_\bullet^{(t)}$ is a graded sequence of ideals on $X$.
\end{definition}

\begin{definition}{\rm (\cite[Definition 4.1]{XZ20b})}
Let $(X,\Delta)$ be an lc log Fano pair and $L=-r(K_X+\Delta)$ a line bundle for some $r\in\bZ^+$. Denote by $R=\oplus_{k\in \bN}R_k:=\oplus_{k\in \bN} H^0(X,kL)$. Suppose $\mF$ is a filtration of $R$, then the \emph{log canonical slope} associated to $\mF$ is defined as follows:
$$\mu(\mF):=\sup \left\{t\in \bR\mid \textit{$\lct\left(X,\Delta; I_\bullet^{(t)}\right)\geq \frac{1}{r}$}\right\}. $$
\end{definition}

\section{Computation of log canonical slope}

In this section, we generalize \cite[Theorem 4.3]{XZ20b} to lc setting.

\begin{theorem}\label{key lemma}
Let $(X,\Delta)$ be an lc log Fano pair and $(\mX,\Delta_\tc;-K_\mX-\Delta_\tc)$ a non-trivial weakly special test configuration. Write $v_{\mX_0}:={\ord_{\mX_0}}|_{K(X)}=c\cdot \ord_E$ for some $c\in \bQ^+$ and some prime divisor $E$ over $X$. Then we have $r\cdot A_{X,\Delta}(E)=\mu(\mF_{\ord_E})$, where $\mF_{\ord_E}$ is an $\bN$-filtration of $R$ defined as follows:
$$\mF^j_{\ord_E}R_k:=\big\{s\in H^0\left(X, -kr\left(K_X+\Delta\right)\right)\mid\textit{$\ord_E(s)\geq j$}\big\}. $$
\end{theorem}

\begin{proof}
We first show $\mu(\mF_{\ord_E})\leq r\cdot A_{X,\Delta}(E)$. Suppose on the contrary that there is a real number $t>r\cdot A_{X,\Delta}(E)$ such that $\lct\left(X,\Delta;I_\bullet^{(t)}\right)\geq \frac{1}{r}$. However, 
$$\lct\left(X,\Delta;I_\bullet^{(t)}\right)\leq \frac{A_{X,\Delta}(E)}{\ord_E\left(I_\bullet^{(t)}\right)}\leq\frac{A_{X,\Delta}(E)}{t}<\frac{1}{r}, $$
which is a contradiction. Note that the second inequality holds since $\ord_E\left(I^{(t)}_\bullet\right)\geq t$. We deal with the converse direction.
By the Claim below, we have 
$$\frac{\mu(\mF_{\ord_E})}{r}-S_{X,\Delta}(E)\geq \fD^{\NA}(\mF_{\ord_E}),$$
see \cite[Definition 2.13]{XZ20b} or \cite[Equation (79)]{Li19} for the definition of $\fD^{\NA}(\mF_{\ord_E})$. It is well-known that 
$$c\cdot \fD^{\NA}(\mF_{\ord_E})=\Ding(\mX,\Delta_\tc; -K_{\mX}-\Delta_\tc),$$
see Remark \ref{rem:Dingeq} where we include a proof. By Remark \ref{coincidence} and Corollary \ref{lc computation}, we see that $\mu(\mF_{\ord_E})\geq r\cdot A_{X,\Delta}(E)$. The proof is finished.
\begin{claim}\label{claim: lcs}
We have the following inequality:
$$\frac{\mu(\mF_{\ord_E})}{r}-S_{X,\Delta}(E)\geq \fD^{\NA}(\mF_{\ord_E}).$$
\end{claim}
\begin{proof}[Proof of the Claim]
Denote by $\mu:=\mu(\mF_{\ord_E})$. We divide the proof into three steps, which deal with the following three cases separately. 
\begin{enumerate}
\item $\mu=0$;
\item $\mu>0$ and $\lct(X,\Delta;I_\bullet^{(t)})>0$ on the interval $(0, \mu+c)$ for some $c>0$;
\item $\mu>0$ and $\lct(X,\Delta;I_\bullet^{(t)})=0$ on the interval $(\mu, \mu+c)$ for any $0<c\ll1$.
\end{enumerate}

\textbf{Step 1.}
We first deal with the case $\mu=0$. In this case, we fix a sufficiently divisible $m$, and first show $\lct(X,\Delta; I_{m,1})=0$. Define $\ka_0=\mO_X$ and $\ka_k:=I_{km,1}$ for $k\in \bZ^+$, it is not hard to see that $\{\ka_k\}_{k\in \bN}$ is a graded sequence of ideals. One may assume that $m$ is sufficiently divisible such that $\ka_k=\ka_{k'}$ for any $k,k'\in \bZ^+$. Observe the following:
\begin{align*}
k\cdot\lct(X,\Delta; \ka_1)=&\ k\cdot\lct(X, \Delta; \ka_k)\\
=&\ \frac{1}{m}\cdot km\cdot \lct(X,\Delta; I_{km, km(\frac{1}{km})})\\
\leq&\  \frac{1}{m}\cdot \lct(X,\Delta; I_\bullet^{(\frac{1}{km})}).
\end{align*}
 If $\lct(X, \Delta; \ka_1)=\lct(X, \Delta; \ka_k)>0$,  one can always choose $k$ sufficiently large such that $\lct(X,\Delta; I_\bullet^{(\frac{1}{km})})>\frac{1}{r}$, which is a contradiction to $\mu=0$.
The contradiction implies that 
$$\lct(X,\Delta; I_{m,1})=\lct(X,\Delta;\ka_1)=0.$$
Then one can choose a valuation $v$ over $X$ computing $\lct(X, \Delta; I_{m,1})$ such that the center of $v$ is contained in the vanishing locus of $I_{m,1}$, thus we see 
$$A_{X,\Delta}(v)=0\quad\text{and}\quad  v(I_{m,1})>0.$$ 
To confirm the inequality in the Claim, we compute 
$$\fD^{\NA}(\mF_{\ord_E})=\fL^{\NA}(\mF_{\ord_E})-\fE^{\NA}(\mF_{\ord_E}),$$ 
see \cite[Equation (75), (78), (79)]{Li19} for $\fD^{\NA}, \fL^{\NA}, \fE^{\NA}$. The term $\fE^{\NA}(\mF_{\ord_E})$ is well-known to be $S_{X,\Delta}(E)$, see \cite[Equation (84)]{Li19}. The term $\fL^{\NA}(\mF_{\ord_E})$ can be expressed as follows by \cite[Equation (72), (73), (78)]{Li19}:
$$\lim_m\ \lct\left(X\times \bA^1, \left(\Delta\times \bA^1\right)\cdot \mI_m^{\frac{1}{mr}}; (t)\right)+\frac{e_+}{r}-1, $$
where 
$$\mI_m:=I_{(m,me_+)}+I_{(m, me_+-1)}t+...+\mO_X\cdot t^{m(e_+-e_-)}\subset I_{m,1}+(t^{m(e_+-e_-)}).$$
By Remark \ref{rem:Dingeq}, it suffices to show the following inequality for the case $\mu=0$ (as we could approximate the above limit at level $m$):
$$0\geq  \lct\left(X\times \bA^1, \left(\Delta\times \bA^1\right)\cdot \mI_m^{\frac{1}{mr}}; (t)\right)+\frac{e_+}{r}-1.$$
We may choose $e_-=0$ in our case, as $\mF_{\ord_E}^0R_k=R_k$ for any $k\in \bN$. Observe the following,
\begin{align*}
0\ =\ &\lct\left(X\times \bA^1, \left(\Delta\times \bA^1\right)\cdot (t^{me_+})^{\frac{1}{mr}}; (t)\right)+\frac{e_+}{r}-1\\
\leq \ &\lct\left(X\times \bA^1, \left(\Delta\times \bA^1\right)\cdot \mI_m^{\frac{1}{mr}}; (t)\right)+\frac{e_+}{r}-1\\
\leq \ &\lct\left(X\times \bA^1, \left(\Delta\times \bA^1\right)\cdot (I_{m,1}+(t^{me_+}))^{\frac{1}{mr}}; (t)\right)+\frac{e_+}{r}-1,
\end{align*}
it suffices to show that
$$\lct\left(X\times \bA^1, \left(\Delta\times \bA^1\right)\cdot (I_{m,1}+(t^{me_+}))^{\frac{1}{mr}}; (t)\right)+\frac{e_+}{r}-1=0.$$
Recall that $v$ is a valuation over $X$ which computes $\lct(X,\Delta; I_{m,1})=0$ with $v(I_{m,1})>0$ and $A_{X,\Delta}(v)=0$. Choose a sufficiently large positive integer $a\gg1$ and define $v_{a, 1}$ to be the quasi-monomial valuation over $X\times \bA^1$ with weights $(a,1)$ along $v$ and $t$. It is not hard to check that $v_{a,1}$ is an lc place of the pair
$$\left(X\times \bA^1, \left(\Delta\times \bA^1\right)\cdot \left(I_{m,1}+(t^{me_+})\right)^{\frac{1}{mr}}\cdot t^{1-\frac{e_+}{r}}\right) ,$$
since 
$$A_{X\times \bA^1, \Delta\times \bA^1}(v_{a,1})-v_{a,1}\left(\left(I_{m,1}+(t^{me_+})\right)^{\frac{1}{mr}}\cdot t^{1-\frac{e_+}{r}}\right)=0. $$
Thus
\begin{align*}
\lct\left(X\times \bA^1, \left(\Delta\times \bA^1\right)\cdot \left(I_{m,1}+(t^{me_+})\right)^{\frac{1}{mr}}; (t)\right)+\frac{e_+}{r}-1=0,
\end{align*}
and therefore
$$0=  \lct\left(X\times \bA^1, \left(\Delta\times \bA^1\right)\cdot \mI_m^{\frac{1}{mr}}; (t)\right)+\frac{e_+}{r}-1. $$
By now, the inequality (actually an equality) in the Claim for the case $\mu=0$ is proved. 

\textbf{Step 2.} From now on, we assume $\mu>0$.
For any real number $b>0$ satisfying $\lct(X,\Delta; I_\bullet^{(t)})>0$ on the interval $t\in (0,b)$,  we first show that $t\mapsto \lct\left(X,\Delta; I_\bullet^{(t)}\right)$ is a non-increasing continuous function on the interval $(0,b)$. This can be proved similarly as the klt setting, as we only deal with positive log canonical thresholds. In the klt setting, the continuity is obtained via the following convexity for $t\in [0,1]$:
$$\frac{1}{\lct(I^tJ^{1-t})}\leq \frac{t}{\lct(I)}+\frac{1-t}{\lct(J)}. $$
In the lc setting, the above relation still holds if we only consider $\lct(I)>0$ and $\lct(J)>0$ (and $\lct(I^tJ^{1-t})>0$ automatically). Let $\mu_0,\mu_1\in (0, b)$ be two rational numbers with $\mu_0<\mu_1$, and $\mu_2=t\mu_0+(1-t)\mu_1$ for some $t\in [0,1]$. Let $m$ be a sufficiently divisible positive integer such that 
$$\lct(X,\Delta; I_\bullet^{(\mu_0)})=m\cdot \lct(X,\Delta; I_{m,m\mu_0}),\ 
\lct(X,\Delta; I_\bullet^{(\mu_1)})=m\cdot \lct(X,\Delta; I_{m,m\mu_1}), \  \text{and}$$
$$\lct(X,\Delta; (I_\bullet^{(\mu_0)})^t\cdot (I_\bullet^{(\mu_1)})^{1-t})=m\cdot\lct(X,\Delta; I_{m,m\mu_0}^t\cdot I_{m,m\mu_1}^{1-t}), $$
where $(I_\bullet^{(\mu_0)})^t$, $(I_\bullet^{(\mu_1)})^{1-t}$ and $(I_\bullet^{(\mu_0)})^t\cdot (I_\bullet^{(\mu_1)})^{1-t}$ are graded sequence of ideals with 
$$(I_\bullet^{(\mu_0)})^t_m=I_{m,m\mu_0}^t, \ (I_\bullet^{(\mu_1)})^t_m=I_{m,m\mu_1}^{1-t},\quad \text{and}\quad ((I_\bullet^{(\mu_0)})^t\cdot (I_\bullet^{(\mu_1)})^{1-t})_m=I_{m,m\mu_0}^t\cdot I_{m,m\mu_1}^{1-t}.$$ 
It is not hard to see the following:
\begin{align*}
\frac{1}{\lct(X,\Delta; I_\bullet^{(\mu_2)})}&\ \leq \frac{1}{\lct(X,\Delta; (I_\bullet^{(\mu_0)})^t\cdot (I_\bullet^{(\mu_1)})^{1-t})}\\
&=\ \frac{1}{m\cdot\lct(X,\Delta; I_{m,m\mu_0}^t\cdot I_{m,m\mu_1}^{1-t})}\\
&\ \leq \frac{t}{m\cdot\lct(X,\Delta; I_{m,m\mu_0})}+ \frac{1-t}{m\cdot \lct(X,\Delta; I_{m,m\mu_1})}\\
&\ = \frac{t}{\lct(X,\Delta; I_\bullet^{(\mu_0)})}+ \frac{1-t}{\lct(X,\Delta; I_\bullet^{(\mu_1)})}.
\end{align*}
This convexity implies the continuity of $t\mapsto \lct\left(X,\Delta; I_\bullet^{(t)}\right)$ on the interval $[\mu_0,\mu_1]$, and hence on $(0,b)$.

Suppose $\lct(X,\Delta; I_\bullet^{(t)})>0$ on the interval $t\in (0, \mu+c)$ for some positive number $c>0$, then the function $t\mapsto \lct\left(X,\Delta; I_\bullet^{(t)}\right)$ is continuous on the interval $(0,\mu+c)$.
By the continuity and the definition of $\mu$, we see the following:
$$\lct\left(X,\Delta; I_\bullet^{(\mu)}
\right)=\frac{1}{r}.$$
Denote by $\ka_m:= I_{m,m\mu}$, then $\lct(X,\Delta; \ka_m)>0$ for a sufficiently divisible $m$.  Consider the graded sequence of ideals $\{\ka_{pm}\}_{p\in \bN}$, then 
$$\lct\left(X,\Delta; I_\bullet^{(\mu)}
\right)=m\cdot \lct(X,\Delta; \ka_{\bullet m})=\frac{1}{r},$$
where $(\ka_{\bullet m})_p=\ka_{pm}$.
Since $\ka_m^p\subset \ka_{pm}$ for any $p\in \bN^+$, we see that the co-support of $\ka_{pm}$ is contained in the co-support of $\ka_m$ for any $p\in \bN^+$. Thus it is enough to consider the valuation whose center is contained in the co-support of $\ka_m$ to compute $\lct(X,\Delta; \ka_{\bullet m})$.
By \cite[Theorem A and Theorem 7.3]{JM12}, one can find a valuation $v$ over $X$ whose center is contained in the co-support of $\ka_m$ such that 
$$A_{X,\Delta}(v)=\frac{1}{r}\cdot v(I_\bullet^{(\mu)}).$$
Since $\lct(X,\Delta; \ka_m)>0$, we see 
$$v(\ka_m)>0 \quad \text{and} \quad \frac{A_{X,\Delta}(v)}{v(\ka_m)}>0. $$
In particular,  $a:=A_{X,\Delta}(v)>0$. By defining $f(\lambda):=v\left(I_\bullet^{(\lambda)}\right)$, we then follow the proof of \cite[Theorem 4.3]{XZ20b} here to conclude this case. First we have $f(\mu)=ar>0$. By the convexity of $f$, we have the following
$$f(\lambda)\geq f(\mu)+\xi(\lambda-\mu)\geq ar+\xi(\lambda-\mu),  $$
where 
$$\xi:=\lim_{h\to 0+}\frac{f(\mu)-f(\mu-h)}{h}. $$
We claim that $\xi>0$. First we have $\xi\geq 0$ since $f(\lambda)$ is non-decreasing. If $\xi=0$, then by the convexity of $f$ we see that $f$ is constant on $(-\infty, \mu]$, which is a contradiction since $f(\mu)=ar>0$ and $f(e_-)=0$. Hence $\xi>0$. Replace $v$ with $\xi^{-1}v$, we may assume $\xi=1$ and then
$$f(\lambda)\geq ar+\lambda-\mu. $$
Let $v_{1,1}$ be the valuation over $X\times \bA^1$ which is the quasi-monomial combination of $v$ and $t$ with weights $(1,1)$, and denote by $e=e_+-e_-$,  then for any $i\in \bN$ we have
\begin{align*}
v_{1,1}(I_{m, me_-+i}\cdot t^{me-i})\ &\geq\  mf(\frac{me_-+i}{m})+(me-i)\\
&\geq\  m(ar+\frac{me_-+i}{m}-\mu)+(me-i)\\
&= \ m(e_++ar-\mu).
\end{align*}
It follows that $v_{1,1}(\mI_m)\geq m(e_++ar-\mu)$ and thus
\begin{align*}
\lct(X\times \bA^1, (\Delta\times \bA^1)\cdot\mI_m^{\frac{1}{mr}}; t)\ &\leq\ A_{X\times \bA^1, \Delta\times \bA^1}(v_{1,1})-\frac{1}{mr}v_{1,1}(\mI_m) \\
&\leq\  a+1-\frac{e_++ar-\mu}{r}\\
&=\ 1-\frac{e_+}{r}+\frac{\mu}{r}.
\end{align*}
Therefore we see
\begin{align*}
\fD^{\NA}(\mF_{\ord_E})&=\fL^{\NA}(\mF_{\ord_E})-\fE^{\NA}(\mF_{\ord_E})\\
&=\lim_m\ \lct\left(X\times \bA^1, \left(\Delta\times \bA^1\right)\cdot \mI_m^{\frac{1}{mr}}; t\right)+\frac{e_+}{r}-1-S_{X,\Delta}(E)\\
&\leq\ \frac{\mu}{r}-S_{X,\Delta}(E).
\end{align*}

\textbf{Step 3}. In this step, we assume $\lct(X,\Delta; I_\bullet^{(t)})=0$ on the interval $t\in (\mu, \mu+c)$ for any $0<c\ll 1$. Let $\mu^+$ be a rational number satisfying $0<\mu^+-\mu\ll1$. Then $\lct(X,\Delta; I_\bullet^{(\mu^+)})=0$. 
By the finite generation of $\mF_{\ord_E}$, one may choose a sufficiently divisible integer $m$ such that 
$$\lct\left(X,\Delta; I_\bullet^{(\mu^+)}\right)=m\cdot \lct(X,\Delta; I_{m,m\mu^+})=0. $$
Therefore one can find a valuation $v$ whose center is contained in the vanishing locus of $I_{m,m\mu^+}$ such that 
$$A_{X,\Delta}(v)=0\quad\text{and}\quad v(I_\bullet^{(\mu^+)})=\frac{v(I_{m,m\mu^+})}{m}>0.$$
Let $v_{a,1}$ be the quasi-monomial combination of $v$ and $t$ with weights $(a,1)$, where $a\gg 1$ is a sufficiently large positive integer, we then have the following inequality:
\begin{align*}
v_{a,1}(\mI_m^{\frac{1}{mr}})&=\ \frac{1}{mr}\cdot v_{a,1}(I_{(m,me_+)}+I_{(m, me_+-1)}t+...+\mO_X\cdot t^{m(e_+-e_-)}) \\
&\geq\ \frac{me_+-m\mu^+}{mr}=\frac{e_+-\mu^+}{r}.
\end{align*}
The inequality holds since 
$$v_{a,1}(I_{m,j}t^{me_+-j})\geq a\cdot v(I_{m,m\mu^+})+me_+-j\gg 1$$ for $me_+\geq j\geq m\mu^+$, and
$$v_{a,1}(I_{m,j}t^{me_+-j})\geq me_+-m\mu^+ $$
for $me_-\leq j\leq m\mu^+$.
Thus we have the following computation:
\begin{align*}
\lct(X\times \bA^1, (\Delta\times \bA^1)\cdot\mI_m^{\frac{1}{mr}}; t)&\leq \ A_{X\times \bA^1, \Delta\times \bA^1}(v_{a,1})-\frac{1}{mr}\cdot v_{a,1}(\mI_m)\\
&\ \leq 1-\frac{e_+-\mu^+}{r}.
\end{align*}
Finally we see 
\begin{align*}
\fD^{\NA}(\mF_{\ord_E})&=\fL^{\NA}(\mF_{\ord_E})-\fE^{\NA}(\mF_{\ord_E})\\
&=\lim_m\ \lct\left(X\times \bA^1, \left(\Delta\times \bA^1\right)\cdot \mI_m^{\frac{1}{mr}}; t\right)+\frac{e_+}{r}-1-S_{X,\Delta}(E)\\
&\leq\ \frac{\mu^+}{r}-S_{X,\Delta}(E).
\end{align*}
As we can choose $\mu^+$ arbitrarily close to $\mu$, thus
$$\fD^{\NA}(\mF_{\ord_E})\leq  \frac{\mu}{r}-S_{X,\Delta}(E).$$
The proof is finished.
\end{proof}
\end{proof}

\begin{remark}\label{rem:Dingeq}
As we cannot find a direct reference which contains the formula
$$c\cdot \fD^{\NA}(\mF_{\ord_E})= \Ding(\mX,\Delta_\tc;\mL:=-K_\mX-\Delta_\tc),$$
we provide a proof here. To compute $\fD^{\NA}(\mF_{\ord_E})$, 
we first construct a sequence of test configurations $\{(\mX_m, \Delta_{\tc,m};\mL_m)\}_m$ for sufficiently divisible $m$.
Define
$$\mI_m:=I_{(m,me_+)}+I_{(m, me_+-1)}t+\cdots+\mO_X\cdot t^{m(e_+-e_-)}, \text{ and}$$
$$I_{(m,x)}:=\Image\left(\mF^x_{v_{\mX_0}}R_m \otimes\mO_X(mr(K_X+\Delta))\to \mO_X\right). $$
Let $\mX_m$ be the normalized blowup of $\mI_m$ on $X\times \bA^1$, and $\mO_{\mX_m}(-E_m)=\mI_m\cdot \mO_{\mX_m}$. Denote by $\pi_m: \mX_m\to X\times \bA^1$ the induced morphism. Define 
$$\mL_m:=-\pi_m^*(K_{X\times \bA^1}+\Delta\times \bA^1)-\frac{1}{mr}E_m,$$ 
then we get a triple $(\mX_m, \Delta_{\tc, m}; \mL_m)\to \bA^1$, where $\Delta_{\tc, m}$ is the natural extension of $\Delta$. By \cite[Lemma 4.6]{Fuj18}, $\mL_m$ is semi-ample for sufficiently divisible $m$. In this case, $(\mX_m, \Delta_{\tc, m}; \mL_m)$ is just a semi test configuration in the sense of \cite[Definition 2.7]{Fuj18}, but we still call it test configuration as we can define the generalized Futaki invariant and the Ding invariant similarly. Note that 
$$H^0(\mX_m, kmr\mL_m)\cong H^0\left(X\times \bA^1, -kmrK_{X\times \bA^1}\otimes \mI_m^k\right)\cong  t^{mke_+}\cdot \bigoplus_{ j\leq mke_+}\left(\mF_{v_{\mX_0}}^jR_{mk}\right)\cdot t^{-j}$$
for $k\in \bN$ and sufficiently divisible $m$ (recall that the filtration $\mF_{v_{\mX_0}}$ has finite generation property).
Thus
\begin{align*}
\Proj \bigoplus_{k\in \bN} H^0(\mX_m, kmr\mL_m)
\cong &\Proj \bigoplus_{k\in \bN}  t^{mke_+}\cdot \bigoplus_{ j\leq mke_+}\left(\mF_{v_{\mX_0}}^jR_{mk}\right)\cdot t^{-j}\\
\cong & \Proj \bigoplus_{k\in \bN}  \bigoplus_{j\in \bZ}\left(\mF_{v_{\mX_0}}^jR_{mk}\right)\cdot t^{-j},
\end{align*}
where one can naturally extend the range of $j$ to $\bZ$. The above relation implies that $(\mX,\Delta_{\tc}; \mL)$ is the ample model of $(\mX_m, \Delta_{\tc,m}; \mL_m)$, denoted by the following diagram:
\[
 \begin{tikzcd}[row sep= 1.25 em]
   \mX_m\arrow[rr,swap,"\phi_m"] \arrow[rd,swap,""]& & \mX\arrow[ld,swap,""]\\
   &\bA^1&
        \end{tikzcd}
\]
Recall the definition of Ding invariant:
$$\Ding(\mX,\Delta_\tc;\mL)=-\frac{\bar{\mL}^{d+1}}{(d+1)(-K_X-\Delta)^d} +\lct\left(\mX, \Delta_\tc+\mD_{\left(\mX,\Delta_\tc;\mL\right)};\mX_0\right)-1,\text{ and} $$
$$\Ding(\mX_m,\Delta_{\tc,m};\mL_m)=-\frac{\bar{\mL}_m^{d+1}}{(d+1)(-K_X-\Delta)^d} +\lct\left(\mX_m, \Delta_{\tc,m}+\mD_{(\mX_m,\Delta_{\tc,m};\mL_m)};\mX_{m,0}\right)-1.$$
Since $\mL_m=\phi_m^* \mL $, we immediately have
$$K_\mX+\mD_{(\mX, \Delta_\tc; \mL)}=\phi_m^*\left(K_{\mX_m}+\mD_{\left(\mX_m,\Delta_{\tc,m};\mL_m\right)}\right).$$
Thus 
$$\lct(\mX, \Delta_\tc+\mD_{(\mX,\Delta_\tc;\mL)};\mX_0)= \lct\left(\mX_m, \Delta_{\tc,m}+\mD_{(\mX_m,\Delta_{\tc,m};\mL_m)};\mX_{m,0}\right),$$
and 
$$\Ding(\mX,\Delta_\tc;\mL)= \Ding(\mX_m,\Delta_{\tc,m};\mL_m).$$
By Fujita's approximation on Ding invariants (see \cite[Equation (91)]{Li19}), we see that
$$\fD^{\NA}\left(\mF_{v_{\mX_0}}\right)=\lim_m \Ding(\mX_m,\Delta_{\tc,m};\mL_m),  $$
which follows that
$$c\cdot \fD^{\NA}(\mF_{\ord_E})= \Ding(\mX,\Delta_\tc;\mL).$$
\end{remark}

\section{Computation of the Generalized Futaki invariants}\label{sec: 3}

In this section, we give a generalization of Theorem \ref{futaki}. The method used here is essentially the same as \cite[Section 5]{Fuj19}, see also \cite{DL20} for a formula of similar type.

\begin{theorem}\label{generalization}
Let $(X,\Delta)$ be a log pair of dimension $n$ and $L$ an ample $\bQ$-line bundle on $X$. Suppose $(\mX,\Delta_\tc;\mL)$ is a non-trivial test configuration of $(X,\Delta;L)$ with integral central fiber. 
Write $\mu(L)=\frac{(-K_X-\Delta).L^{n-1}}{L^n}$ and $v_{\mX_0}:={\ord_{\mX_0}}|_{K(X)}=c\cdot \ord_E$ for some prime divisor $E$ over $X$ and $c\in \bQ^+$, then we have the following reformulation:
$$\Fut(\mX,\Delta_\tc;\mL)=A_{X,\Delta}(v_{\mX_0})-\mu(L)S_{X,L}(v_{\mX_0}) + P(\mX,\Delta_\tc;\mL).$$
Here, $A_{X,\Delta}(v_{\mX_0})=c\cdot A_{X,\Delta}(E)$, $S_{X,L}(v_{\mX_0})=c\cdot S_{X,L}(E)$, and 
$$P(\mX,\Delta_\tc;\mL)= \frac{1}{L^n}\left(\tau^*\bar{\mL}\right)^n\cdot\rho^*\left(K_{X\times \bP^1/\bP^1}+\Delta\times \bP^1+\mu\left(L\right)L_{X\times \bP^1}\right),$$
where $L_{X\times \bP^1}:=L\times \bP^1$,
$\rho$ and $\tau$ come from the following normalization of the graph:
\[
 \begin{tikzcd}[row sep= 1.25 em]
   & \mY   \arrow[rd,"\tau"] \arrow[ld,swap,"\rho"]  &  \\
   X \times \bP^1 \arrow[rr,dashed]\arrow[rd,swap,"p"] & & \bar{\mX}\arrow[ld,swap]\\
   & \bP^1 &
        \end{tikzcd}
        \]        
\end{theorem}

\begin{proof}
We first write down the following formulas:
$$K_{\mY/\bP^1}+\rho_*^{-1}\left(\Delta\times\bP^1\right)=\rho^*\left(K_{X\times \bP^1/\bP^1}+\Delta\times \bP^1\right)+b\tilde{\mX_0}+F,\text{ and} $$
$$K_{\mY/\bP^1}+\tau_*^{-1}\bar{\Delta}_\tc=\tau^*\left(K_{\bar{\mX}/\bP^1}+\bar{\Delta}_\tc\right)+a \tilde{X}+G,  $$
where $\tilde{\mX_0}$ (resp. $\tilde{X}$) is the strict transformation of $\mX_0$ (resp. $X\times 0$) on $\mY$, and $F$ and $G$ contain components in the central fiber of $\mY\to \bP^1$ other than $\tilde{\mX}_0$ and $\tilde{X}$.
Note that $F$  and $G$ are both $\tau$-exceptional and $\rho$-exceptional. It is clear that $b=A_{X,\Delta}(v_{\mX_0})$ (see \cite[Proposition 4.11]{BHJ17}). 
We also write
\begin{align}\label{difference}
\tau^*\bar{\mL}-\rho^*L_{X\times \bP^1} =a_1 \tilde{X}+b_1 \tilde{\mX_0}+Q,
\end{align}
where $Q$ is a divisor supported in the central fiber of $\mY\to \bP^1$ which is both $\tau$-exceptional and $\rho$-exceptional.
Then we can see that
\begin{align}\label{b1}
b_1=\frac{1}{L^n}\left(\tau^*\bar{\mL}\right)^n\cdot\left(\tau^*\bar{\mL}-\rho^*L_{X\times \bP^1}\right)=\frac{1}{L^n}\bar{\mL}^{n+1}-\frac{1}{L^n}\left(\tau^*\bar{\mL}\right)^n\cdot\rho^*L_{X\times \bP^1} .
\end{align}
Recall that
\begin{align*}
\Fut(\mX,\Delta_\tc;\mL)&=\frac{\bar{\mL}^n\cdot \left(K_{\bar{\mX}/\bP^1}+\bar{\Delta}_\tc\right)}{L^n} -\frac{n}{n+1}\cdot\frac{\left(\left(K_X+\Delta\right)L^{n-1}\right)\bar{\mL}^{n+1}}{(L^n)^2}\\
&=\frac{\bar{\mL}^n\cdot \left(K_{\bar{\mX}/\bP^1}+\bar{\Delta}_\tc\right)}{L^n}+\frac{n}{n+1}\mu(L)\frac{\bar{\mL}^{n+1}}{L^n}.
\end{align*}
To reformulate the generalized Futaki invariant $\Fut(\mX,\Delta_\tc;\mL)$, we first calculate $\frac{\bar{\mL}^n\cdot\left(K_{\bar{\mX}/\bP^1}+\bar{\Delta}_\tc\right)}{L^n}$ as follows:
\begin{align}\label{first term}
\frac{\bar{\mL}^n\cdot\left(K_{\bar{\mX}/\bP^1}+\bar{\Delta}_\tc\right)}{L^n}=&\frac{1}{L^n}\left(\tau^*\bar{\mL}\right)^n\cdot\tau^*\left(K_{\bar{\mX}/\bP^1}+\bar{\Delta}_\tc\right)\\
=&\frac{1}{L^n}\left(\tau^*\bar{\mL}\right)^n\cdot \left(\rho^*(K_{X\times \bP^1/\bP^1}+\Delta\times \bP^1)+b\tilde{\mX_0}-a\tX\right)\\
=&\frac{1}{L^n}\left(\tau^*\bar{\mL}\right)^n\cdot \rho^*\left(K_{X\times \bP^1/\bP^1}+\Delta\times\bP^1\right)+A_{X,\Delta}(v_{\mX_0}).
\end{align}
Next we compute $\frac{n}{n+1}\mu(L)\frac{\bar{\mL}^{n+1}}{L^n}$.

\begin{claim}\label{key claim}
Let $D$ be an effective divisor on $X$ corresponding to $f\in H^0(X,mL)$. Then $f\in \mF^j_{v_{\mX_0}}H^0(X,mL)$ if and only if $f\in \mF^{j+mb_1}_{\ord_{\mX_0}}H^0(\mX,m\mL)$, where
$$\mF^j_{v_{\mX_0}}H^0(X,mL):=\big\{f\in H^0(X,mL)\mid\textit{$v_{\mX_0}(f)\geq j$}\big\}, $$
$$\mF^{j}_{\ord_{\mX_0}}H^0(\mX,m\mL):=\left\{f\in H^0(X,mL)\mid\textit{$t^{-j}\cdot \rho^*p_1^*f\in H^0(\mX,m\mL)$}\right\},$$
and $p_1: X\times \bA^1\to X$ is the projection to $X$.
\end{claim}

\begin{proof}[Proof of the Claim]
It is clear that $f\in \mF^{j}_{\ord_{\mX_0}}H^0(\mX,m\mL)$ if and only if the following holds:
$$\rho^*p_1^*D+\tau^*m\mL-\rho^*mL_{X\times \bA^1}-j\tau^*\mX_0\geq 0. $$
As 
$$\rho^*p_1^*D+\tau^*m\mL-\rho^*mL_{X\times \bA^1}-j\tau^*\mX_0\sim_\bQ \tau^*(m\mL-j\mX_0),$$ 
the above inequality is also equivalent to
$$\tau_*\left( \rho^*p_1^*D+\tau^*m\mL-\rho^*mL_{X\times \bA^1}-j\tau^*\mX_0\right) \geq 0.$$
Combine (\ref{difference}), this is equivalent to the condition
$$\ord_{\tilde{\mX}_0}\left(\rho^*p_1^*D\right)\geq j-mb_1. $$
Thus the claim holds.
\end{proof}

Armed by the above claim, we can calculate $w(m)$, which is the total weights of $\bC^*$-action on $H^0(\mX_0,m\mL_0)$:
\begin{align*}
w(m)&=\sum_{j\in \bZ} j \left(\dim \mF^j_{\ord_{\mX_0}}H^0\left(\mX,m\mL\right)-\dim \mF^{j+1}_{\ord_{\mX_0}}H^0\left(\mX,m\mL\right)\right)\\
&=\sum_{j\geq 0} (j+mb_1)\left(\dim \mF^j_{v_{\mX_0}}H^0\left(X,mL\right)-\dim \mF^{j+1}_{v_{\mX_0}}H^0\left(X,mL\right)\right)\\
&=\sum_{j\geq 0} j \left(\dim \mF^j_{v_{\mX_0}}H^0(X,mL)-\dim \mF^{j+1}_{v_{\mX_0}}H^0(X,mL)\right) +mb_1 \dim H^0(X,mL)\\
&=\sum_{j\geq 1}\dim \mF^j_{v_{\mX_0}}H^0(X,mL)+mb_1 \dim H^0(X,mL) .
\end{align*}
Thus we get
$$\frac{\bar{\mL}^{n+1}}{L^n}=\lim_m \frac{1}{L^n}\frac{w(m)}{m^{n+1}/(n+1)!}= (n+1)\big(S_{X,L}(v_{\mX_0})+b_1\big),$$
where $S_{X,L}(v_{\mX_0})=c\cdot S_{X,L}(E)$. Combine the expression of $b_1$ (see (\ref{b1})), we have 
\begin{align}
b_1=\frac{1}{n}\frac{1}{L^n}\cdot \left(\tau^*\bar{\mL}\right)^n\cdot\rho^*L_{X\times \bP^1} -\frac{n+1}{n}S_{X,L}(v_{\mX_0}).
\end{align}
Then
\begin{align}\label{mu-term}
\frac{n}{n+1}\mu(L)\cdot\frac{\bar{\mL}^{n+1}}{L^n}=\frac{\mu(L)}{L^n}\left(\tau^*\bar{\mL}\right)^n\cdot\rho^*L_{X\times \bP^1}-\mu(L)S_{X,L}(v_{\mX_0}).
\end{align}
Combine (\ref{first term}) and (\ref{mu-term}) the proof is finished.
\end{proof}

\section{Dreamy property}
In this section, we prove Theorem \ref{thm: main2}.  Let us begin with the following lemma (see also \cite[Theorem 1]{Moraga20}).

\begin{lemma}\label{extraction}
Let $(X,\Delta)$ be an lc log pair and $E$ an lc place. Then one can find an extraction birational morphism $Y''\to X$ such that $-E$ is ample over $X$. In particular, if the center of $E$ is of at least codimension 2, then $E$ is the unique exceptional divisor.
\end{lemma}

\begin{proof}
Take a dlt modification $f: (Y, \Delta_Y)\to (X, \Delta)$ such that $(Y,\Delta_Y)$ is $\bQ$-factorial dlt, $K_Y+\Delta_Y=f^*(K_X+\Delta)$,  and $E$ is a component of $\Delta_Y^{=1}$ (i.e., the sum of all components of $\Delta_Y$ whose coefficients are one), see \cite[Theorem 3.1]{KK10} or \cite[Theorem 10.5]{Fujino11}. Choose a small rational number $0<\epsilon<1$ and consider the pair $(Y,\Delta_Y-\epsilon E)$. It is clear that $(Y, \Delta_Y-\epsilon E)$ is $\bQ$-factorial dlt and $(K_Y+\Delta_Y-\epsilon E)\sim_{\bQ, f} -\epsilon E$. By \cite[Theorem 1.1]{Birkar12} or \cite[Theorem 1.6]{HX13}, one can run a $(K_Y+\Delta_Y-\epsilon E)$-MMP over $X$ with scaling, which terminates with a good minimal model $Y\dashrightarrow Y'$ over $X$. Let $Y'\to Y''$ be the ample model over $X$ and $E''$ the push-forward of $E$ to $Y''$. We claim that $E$ is not contracted, i.e., $E''\ne 0$. Denote $\Delta_{Y''}$ to be the push-forward of $\Delta_Y$, assume $E$ is contracted, then we have the following inequality for log discrepancies: 
$$\epsilon=A_{Y, \Delta_Y-\epsilon E}(E)<A_{Y'', \Delta_{Y''}}(E'')=A_{Y,\Delta_Y}(E)=0, $$
a contradiction. Thus $E''\ne 0$ and $-E''$ is ample over $X$.  For any curve $C$ on $Y''$ which is contracted by $Y''\to X$, we have $E''.C<0$, thus the locus of $C$ lies inside $E''$. If $E''$ is exceptional, this implies that $E''$ is the unique exceptional divisor of $Y''\to X$. The proof is finished.
\end{proof}

We are ready to prove Theorem \ref{thm: main2}.

\begin{proof}[Proof of Theorem \ref{thm: main2}]
Let $D\sim_\bQ -K_X-\Delta$ be an effective $\bQ$-divisor on $X$ such that the log pair $(X,\Delta+D)$ is lc and $E$ is an lc place of $(X,\Delta+D)$. Choose a divisible positive integer $r$ such that both $L:=-r(K_X+\Delta)$ and $r(K_X+\Delta+D)$ are Cartier. Let $Z$ be the affine cone over $X$ with respect to $L$, i.e., 
$$Z=\Spec \bigoplus_{m\in \bN}H^0(X, -mr(K_X+\Delta)),$$ 
and $o\in Z$ is the cone vertex. Denote $\Delta_Z$ and $D_Z$ to be the extensions of $\Delta$ and $D$ on $Z$ respectively. Consider the projection morphism $p: Z\setminus o\to X$, we denote $E_\infty$ to be the prime divisor over $Z\setminus o$ via pulling back $E$ over $X$. Let $\mu: \tZ\to Z$ be the blowup of $Z$ at the vertex and $X_0$ the exceptional divisor, then we write $v_0:=\ord_{X_0}$ for the canonical valuation. It is clear that $E_\infty$ is an lc place of the pair $(Z, \Delta_Z+D_Z)$, and by \cite[Proposition 3.14]{Kollar13}, $v_0$ is also an lc place of $(Z, \Delta_Z+D_Z)$. Write $w_k:=k\cdot v_0+\ord_{E_\infty}$ to be the quasi-monomial valuation with weights $(k,1)$ alone $v_0$ and $\ord_{E_\infty}$ for some positive integer $k$. By \cite[Proposition 5.1]{JM12}, $w_k$ is an lc place of $(Z,\Delta_Z+D_Z)$ whose center is exactly the cone vertex $o\in Z$. By Lemma \ref{extraction}, there is an extraction $f: W\to Z$ which has a unique exceptional divisor $E_k$ corresponding to $w_k$ and $-E_k$ is $f$-ample. We write
$$K_W+f_*^{-1}\Delta_Z+f_*^{-1}D_Z+E_k=f^*(K_Z+\Delta_Z+D_Z). $$
Consider the following exact sequence for $pl\in \bN$, where $l\in \bN$ and $p$ is a fixed divisible positive integer such that $pE_k$ is Cartier:
$$0\to \mO_{W}(-(pl+1)E_k)\to \mO_{W}(-plE_k)\to \mO_{E_k}(-plE_k)\to 0. $$
Write $\ka_l:=f_*\mO_{W}(-lE_k)$, then 
$$\ka_{pl}/\ka_{pl+1}\cong H^0(E_k, -plE_k|_{E_k})$$ 
for $l\in \bN$, since $-E_k$ is $f$-ample and $R^1 f_*\mO_{W}(-(pl+1)E_k)=0$ by an lc version of Kawamata-Viehweg vanishing (see \cite[Theorem 1.7]{Fujino14}). Thus the graded algebra $\oplus_{l\in \bN}\ka_{pl}/\ka_{pl+1}$ is finitely generated, hence so is $\oplus_{l\in \bN}\ka_l/\ka_{l+1}$.  Therefore, the graded algebra $\oplus_{l\in \bN}\ka_l$ is finitely generated.
Denote by 
$$R_m:=H^0(X,-mr(K_X+\Delta)) \quad \text{and}\quad  R:=\oplus_{m\in \bN}R_m.$$ 
We naturally extend $\oplus_{l\in \bN}\ka_l$ and $R$ to be graded algebras indexed by $\bZ$ via defining $R_m=0$ for $m<0$ and $\ka_l=\mO_Z$ for $l< 0$. We use the notation 
$$\mF^j_{\ord_E}R_m:=\{s\in R_m\mid  \textit{$\ord_E(s)\geq j$}\}\quad \text{and} \quad \ka_{m,l}:=\{s\in R_m\mid \textit{$km+\ord_E(s)\geq l$}\} .$$
It is clear that $\ka_l=\oplus_{m\in \bZ}\ka_{m,l}$, and $s\in \mF^j_{\ord_E}R_m$ if and only if $s\in \ka_{m,km+j}$.
Then we have
$$\bigoplus_{m\in \bZ}\bigoplus_{j\in \bZ} \mF^j_{\ord_E}R_m\cong \oplus_{m\in \bZ}\oplus_{j\in \bZ} \ka_{m,km+j}\cong \oplus_{m\in \bZ}\oplus_{j\in \bZ} \ka_{m,j}\cong \oplus_{j\in \bZ} \ka_j,$$
which are all finitely generated graded algebras.
Note that 
$$\mF^j_{\ord_E}R_m\cong H^0(Y, g^*(-mr(K_X+\Delta))-jE ), $$
where $g: Y\to X$ is a projective normal birational model such that $E$ is a prime divisor on $Y$, thus $\bigoplus_{m\in \bN}\bigoplus_{j\in \bN} \mF^j_{\ord_E}R_m$ being finitely generated implies that $E$ is dreamy. The proof is finished.
\end{proof}

\section{Correspondence}

In this section, we prove Theorem \ref{thm: main1}. We will divide it into two parts, i.e., Theorem \ref{thm: tc} and Theorem \ref{thm: complement}. 

\begin{theorem}\label{thm: tc}
Let $(X,\Delta)$ be a strictly lc log Fano pair. Suppose $E$ is an lc place of complements of $(X,\Delta)$, then $E$ induces a non-trivial weakly special test configuration  $(\mX, \Delta_\tc; -K_{\mX}-\Delta_\tc)$ of $(X,\Delta; -K_X-\Delta)$ such that the restriction of $\ord_{\mX_0}$ to the function field $K(X)$ is  exactly $\ord_E$.
\end{theorem}

\begin{proof}
By assumption, there exists a $\bQ$-divisor $D\ge0$ on $X$ such that $ D\sim_\bQ -K_X-\Delta$ and $E$ is an lc place of $(X,\Delta+D)$. If $E$ is exceptional over $X$, then by Lemma \ref{extraction}, there exists an extraction $f: Y\to X$ which only extracts $E$, and we have
$$K_Y+f_*^{-1}\Delta+f^{-1}_*D+E=f^*(K_X+\Delta+D). $$
If $E$ is a prime divisor on $X$, then we may just take $f:={\rm id}:Y:=X\to X$. 
We only deal with the case where $E$ is exceptional in the remaining proof, as one can prove similarly for the case where $E$ is a divisor on $X$.
According to Theorem \ref{thm: main2}, we know that the following $\bZ^2$-graded ring 
$$\bigoplus_{k\in \bN}\bigoplus_{j\in \bN} H^0(Y, f^*(-kr(K_X+\Delta))-jE) $$
is finitely generated for some $r\in \bZ^+$ such that $-r(K_X+\Delta)$ is Cartier. Let 
$$R_k:=H^0(X, -kr(K_X+\Delta)) \text{ and } \mF^j_{\ord_E}R_k:=H^0(Y, f^*(-kr(K_X+\Delta))-jE).$$ 
We construct the following degeneration family over $\bA^1$:
$$\mZ:=\Proj \bigoplus_{k\in \bN}\bigoplus_{j\in \bZ} \left(\mF^j_{\ord_E}R_k\right) t^{-j}\to \bA^1.$$
 By \cite[Lemma 3.8]{Fuj17} we know that $(\mZ, \Delta_\mZ; -K_\mZ-\Delta_\mZ)$ is a test configuration of $(X,\Delta)$ whose central fiber $\mZ_0$ is integral.  Let $\Delta_\mZ$ and $\mD_\mZ$ be the extensions of $\Delta$ and $D$ on $\mZ$ respectively.  Then it holds that ${\ord_{\mZ_0}}|_{K(X)}= \ord_E$\footnote{By Claim \ref{key claim}, we see that a test configuration $(\mX,\Delta_\tc;\mL)$ (whose central fiber is integral) with $v_{\mX_0}:={\ord_{\mX_0}}|_{K(X)}=c\cdot \ord_E$ for some $c\in \bQ^+$ and prime divisor $E$ over $X$
 is induced by the filtration $\mF_{v_{\mX_0}}$. Since our test configuration here is induced by $\mF_{\ord_E}$, we have $c=1$.}. 
We next show that $(\mZ, \Delta_{\mZ}+\mZ_0)$ is log canonical, which will imply that $(\mZ, \Delta_{\mZ}; -K_{\mZ}-\Delta_\mZ)$ is a weakly special test configuration of $(X,\Delta)$ induced by $E$.

Put the morphism $f$ into the trivial family $F: Y_{\bA^1}\to X_{\bA^1}$, then we have
$$K_{Y_{\bA^1}}+F_*^{-1}\Delta_{\bA^1}+F_*^{-1}D_{\bA^1}+E_{\bA^1}+Y_0=F^*\left(K_{X_{\bA^1}}+\Delta_{\bA^1}+D_{\bA^1}+X_0\right) ,$$
where $X_0$ (resp. $Y_0$) is the central fiber of the family $X_{\bA^1}\to \bA^1$ (resp. $Y_{\bA^1}\to \bA^1$). 
Let $v$ be the divisorial valuation over $X_{\bA^1}$ with weights $(1,1)$ along divisors $E_{\bA^1}$ and $Y_0$, then 
$$A_{X_{\bA^1}, \Delta_{\bA^1}+D_{\bA^1}+X_0}(v)=A_{Y_{\bA^1}, F_*^{-1}\Delta_{\bA^1}+F_*^{-1}D_{\bA^1}+E_{\bA^1}+Y_0}(v)=0. $$
By Lemma \ref{extraction}, one can extract $v$ to be a divisor on a birational model $g: \mY\to X_{\bA^1}$ as follows:
$$K_{\mY}+g_*^{-1}\Delta_{\bA^1}+g_*^{-1}X_0+g_*^{-1}D_{\bA^1}+\mE=g^*(K_{X_{\bA^1}}+\Delta_{\bA^1}+D_{\bA^1}+X_0),$$
where $v=\ord_\mE$.
We note here that the pair $\left(\mY, g_*^{-1}\Delta_{\bA^1}+g_*^{-1}D_{\bA^1}+g_*^{-1}X_0+\mE\right)$ is  log canonical and log Calabi-Yau over $\bA^1$. 

Now we may compare the two pairs, 
$$\left(\mY, g_*^{-1}\Delta_{\bA^1}+g_*^{-1}D_{\bA^1}+g_*^{-1}X_0+\mE\right) \quad \text{and}\quad (\mZ, \Delta_{\mZ}+D_{\mZ}+\mZ_0).$$ 
Since $\mE$ and $\mZ_0$ induce the same divisorial valuation on $K(\mY)\cong K\left(X\times \bA^1\right)\cong K(\mZ)$, we have the following birational contraction map:
$$\left(\mY, g_*^{-1}\Delta_{\bA^1}+g_*^{-1}D_{\bA^1}+g_*^{-1}X_0+\mE\right)\dashrightarrow (\mZ, \Delta_{\mZ}+D_{\mZ}+\mZ_0). $$
As both two sides are log Calabi-Yau over $\bA^1$, they are crepant\footnote{Note here that $\mY\dashrightarrow \mZ$ is a birational contraction, thus the crepant property follows from non-negativity lemma, e.g., \cite[Lemma 3.39]{KM98}. Here one might argue more that $K_\mZ+\Delta_\mZ+D_\mZ$ is $\bQ$-Cartier. As we have mentioned,  $(\mZ,\Delta_{\mZ})$ is a test configuration with integral central fiber, thus we naturally have that both $D_\mZ$ and $K_{\mZ}+\Delta_\mZ$ are $\bQ$-Cartier. }, i.e., there is a common log resolution $p: \mW\to \mY$ and $q: \mW\to \mZ$ such that
$$p^*\left(K_{\mY}+g_*^{-1}\Delta_{\bA^1}+g_*^{-1}D_{\bA^1}+g_*^{-1}X_0+\mE\right)=q^*(K_{\mZ}+\Delta_{\mZ}+D_{\mZ}+\mZ_0). $$
It follows that $(\mZ, \Delta_{\mZ}+D_{\mZ}+\mZ_0)$ is log canonical, and thus so is $(\mZ, \Delta_{\mZ}+\mZ_0)$. This completes the proof.
\end{proof}

\begin{theorem}\label{thm: complement}
Let $(X,\Delta)$ be a strictly lc log Fano pair.  Suppose $(\mX,\Delta_\tc; -K_{\mX}-\Delta_\tc)$ is a non-trivial weakly special test configuration of $(X,\Delta;-K_X-\Delta)$, then the restriction of $\ord_{\mX_0}$ to $K(X)$ is an lc place of complements of $(X,\Delta)$.
\end{theorem}

\begin{proof}
Armed by Theorem \ref{key lemma}, we can argue in the same way as the proof of \cite[Theorem 4.10]{Xu21} for this direction. We use the notation of Theorem \ref{key lemma}. We may assume $A_{X,\Delta}(E)>0$\footnote{There is nothing to prove for the case when $A_{X,\Delta}(E)=0$, as then $E$ is already an lc place of $(X,\Delta)$.}. We first show that the function $t\mapsto \lct\left(X,\Delta;I_\bullet^{(t)}\right)$ is continuous on the interval $(0,\mu]$, where $\mu:=\mu(\mF_{\ord_E})=r\cdot A_{X,\Delta}(E)$. By the same discussion as in Step 2 of Claim \ref{claim: lcs}, it is enough to show $\lct\left(X,\Delta;I_\bullet^{(\mu)}\right)>0$. Note here that $\mu>0$ is a rational number.
Suppose $\lct\left(X,\Delta;I_\bullet^{(\mu)}\right)=0$, 
we have the following relation for sufficiently divisible integer $m$: 
$$\lct\left(X,\Delta; I_\bullet^{(\mu)}\right)=m\cdot \lct(X,\Delta; I_{m,m\mu})=0. $$
Therefore one can find a valuation $v$ whose center is contained in the vanishing locus of $I_{m,m\mu}$ such that 
$$A_{X,\Delta}(v)=0 \quad \text{and}\quad \frac{1}{m}v(I_{m,m\mu})=v(I_\bullet^{(\mu)})>0.$$
By the continuity of $t\mapsto v(I_\bullet^{(t)})$, we see that $v(I_\bullet^{(\mu')})>0$ for some rational number $\mu'$ with $0<\mu-\mu'\ll1$. However, $A_{X,\Delta}(v)=0$ together with $v(I_\bullet^{(\mu')})>0$
imply that $\lct (X,\Delta; I_\bullet^{(\mu')})=0$, which is a contradiction to the definition of $\mu$. Thus the function $t\mapsto \lct\left(X,\Delta;I_\bullet^{(t)}\right)$ is continuous on the interval $(0,\mu]$.
By the continuity we know 
$$\lct\left(X,\Delta;I_\bullet^{(r\cdot A_{X,\Delta}(E))}\right)\geq \frac{1}{r}. $$
However, 
$$\lct\left(X,\Delta;I_\bullet^{(r\cdot A_{X,\Delta}(E))}\right)\leq \frac{A_{X,\Delta}(E)}{\ord_E\left(I_\bullet^{(r\cdot A_{X,\Delta}(E))}\right)} =\frac{1}{r}.$$
Thus we have
$$\lct\left(X,\Delta;I_\bullet^{(r\cdot A_{X,\Delta}(E))}\right)= \frac{1}{r}. $$
As the filtration is finitely generated (see \cite[Proposition 2.15]{BHJ17}), one can choose a sufficiently divisible $k\in \bZ^+$ such that
$$\lct\left(X,\Delta;I_\bullet^{(r\cdot A_{X,\Delta}(E))}\right)=k\cdot\lct\left(X,\Delta;I_{k,krA_{X,\Delta}(E)}\right)=\frac{1}{r}. $$
This means that there is a $D\in |-kr(K_X+\Delta)|$ such that $\left(X,\Delta+\frac{1}{kr}D\right)$ is lc and $\ord_E(D)\geq kr\cdot A_{X,\Delta}(E)$, which implies that $E$ is an lc place of complements. 
\end{proof}

\begin{proof}[Proof of Theorem \ref{thm: main1}]
The proof is a combination of Theorem \ref{thm: tc} and Theorem \ref{thm: complement}.
\end{proof}

\begin{corollary}
Let $(X,\Delta)$ be a strictly lc log Fano pair and denote by $R=\oplus_{k\in \bN}R_k:=\oplus_{k\in \bN}H^0(X,-kr(K_X+\Delta))$. A prime divisor $E$ over $X$ is an lc place of complements of $(X,\Delta)$ if and only if $E$ is dreamy and $\mu(\mF_{\ord_E})=r\cdot A_{X,\Delta}(E)$. 
\end{corollary}

\begin{proof}
The proof is the same as that of \cite[Lemma A.7]{XZ20b} just by applying Theorem \ref{thm: main1} and Theorem \ref{thm: main2} at the corresponding places.
\end{proof}

\bibliographystyle{amsalpha}
\bibliography{reference.bib}
\end{document}